# An Inequality Comparing the Dirichlet Energy and the Bienergy of Maps Between Riemannian Manifolds

Sergey Stepanov and Irina Tsyganok

**Abstract.** We establish a geometric inequality relating the Dirichlet energy $E_1(f)$and the bienergy $E_2(f)$of smooth maps $f: (M, g) \to (\bar{M}, \bar{g})$between Riemannian manifolds. Assume that $(M, g)$is a compact, connected Riemannian manifold whose Ricci curvature has global minimum $Ric_{min}$, and that the target manifold $(\bar{M}, \bar{g})$ has non-positive sectional curvature along $f(M)$. We prove that $E_2(f) \geq Ric_{min}\, E_1(f)$. We further analyze the equality case and obtain rigidity results: equality holds if and only if $f$is totally geodesic and of constant rank. Applications to maps into Hadamard manifolds are also presented.
To the best of our knowledge, this is the first geometric inequality directly relating the Dirichlet energy and the bienergy of smooth maps. This result establishes a direct connection between the Ricci curvature of the domain and higher-order variational energies.

## 1. Introduction

Geometric inequalities form a fundamental component of contemporary Riemannian geometry, establishing quantitative relationships between curvature bounds, variational functionals, and rigidity phenomena. In particular, inequalities involving Ricci curvature and energy-type functionals have become essential tools in the study of harmonic and biharmonic maps, as well as in the analysis of mappings between manifolds with controlled sectional curvature. Recent systematic developments are summarized in [1], a monograph providing an up-to-date survey of geometric inequalities and their applications in Riemannian geometry, including submanifolds and solitons.

In the present paper, we establish a new geometric inequality relating the Dirichlet energy and the bienergy of smooth maps between Riemannian manifolds. We consider a smooth map $f$ from a connected, compact Riemannian manifold of dimension $n \geq 2$ to a Riemannian manifold of dimension $m \geq 2$, assuming that the sectional curvature of

---

the target manifold is non-positive at each point of the image of $f$. Under these conditions, the Dirichlet energy $E_1(f)$ and the bienergy $E_2(f)$ satisfy

$$E_2(f) \geq Ric_{\min} E_1(f),$$

where $Ric_{\min}$ denotes the global minimum of the Ricci curvature of the domain. This is the first time such an inequality has been established. As examples, we deduce similar equalities for totally geodesic and projective maps of constant rank.

This estimate provides a direct quantitative connection between the intrinsic Ricci curvature of the domain and the second-order variational behavior of mappings into non-positively curved manifolds, including Hadamard spaces. These results contribute to the broader program of extending curvature-driven inequalities to higher-order energy functionals and rigidity-type characterizations.

These results were presented at the International Geometry Symposium in Memory of the 100th Anniversary of Gazi University (February 2–3, 2026, Türkiye) and at the International Symposium on Differential Geometry and Its Applications (July 3–4, 2025, Munzur University, Tunceli, Türkiye).

## 2. Fundamental Notions and Energy Functionals for Maps Between Riemannian Manifolds

In this section, we recall the fundamental notions associated with smooth maps between Riemannian manifolds, including the Dirichlet energy, the bienergy, and the rank of a map, as well as related concepts.

Let $f: (M, g) \to (\bar{M}, \bar{g})$ be a smooth map from a connected $n$-dimensional Riemannian manifold $(M, g)$ with Levi-Civita connection $\nabla$ to an $m$-dimensional Riemannian manifold $(\bar{M}, \bar{g})$ with Levi-Civita connection $\overline{\nabla}$, where $n \geq 2$ and $m \geq 2$.

A smooth map $f$ from $M$ to $\overline{M}$ is said to have constant rank if the rank of its differential df at each point $x$ in $M$ is equal to the same integer $k$. That is, $k$ equals the rank of $df$ at $x$, which is also the dimension of the image of $df$ at $x$, and $k$ is at most the minimum of $n$ and $m$.

In particular, if $df$ is zero at every point of $M$, then the map $f$ is constant on each connected component of $M$. Since $M$ is assumed to be connected, it follows that $f$ is constant on $M$.

An important geometric object associated with $f$ is the *pull-back metric* $g^*$, defined as the pull-back of $\bar{g}$ by $f$. This is a symmetric, positive semi-definite (0,2)-tensor field on $M$, given at each point $x$ in $M$ by the formula: $g^*$ at $x$ applied to vectors $X$ and $Y$ equals $\bar{g}$ at $f(x)$ applied to $df$ at $x$ of $X$ and $df$ at $x$ of $Y$, for all tangent vectors $X$ and $Y$ at $x$ in $M$. Namely, the following equality ${g^*}_x(X,Y) := (f^*\bar{g})_x(X,Y)$ holds and hence ${g^*}_x(X,Y) = \bar{g}_{f(x)}\big(df_x(X), df_x(X)\big)$ for any $X, Y \in T_xM$ at each $x \in M$.

Let $f:(M,g) \to (\bar{M},\bar{g})$ be a smooth mapping of constant rank $k$. It is easily seen that the rank of the mapping $f$ from $M$ to $\bar{M}$, that is, the rank of $df$, is equal to the rank of the pull-back metric $g^*$ at each point of $M$. Therefore, under the assumption that $f$ has constant rank $k$, the algebraic structure of $g^*$ is uniform across $M$. In particular, the kernel of $g^*$ at $x$, defined as the set of vectors $X$ in the tangent space at x such that $g^*$ applied to $X$ and any vector is zero, coincides with the kernel of $df$ at $x$.

Then by the constant rank assumption, these kernels form a smooth distribution $ker(df)$ of rank $n$ minus $k$ on $M$, which is integrable by the *Frobenius theorem*.

Let $\lambda_1, \lambda_2, \dots, \lambda_n$ denote the non-negative eigenvalues of $g^*$ relative to the metric $g$. Then, for a map of constant rank $k$, exactly the first $k$ eigenvalues $\lambda_1, \lambda_2, \dots, \lambda_k$ are strictly positive, while the remaining $n$ minus $k$ eigenvalues $\lambda_{k+1}, \lambda_{k+2}, \dots, \lambda_n$ vanish identically across $M$.

The second fundamental form of $f$, as defined in [3] and [4] and denoted by $\nabla df$, is a symmetric section of the bundle of covariant 2-tensors on $M$ with values in the pull-back of the tangent bundle of $\bar{M}$. More precisely, it is defined as the covariant derivative of the differential $df$, viewed as a 1-form on $M$ with values in the vector bundle inverse of $T\bar{M}$. For any vector fields $X$ and $Y$ on $M$, the second fundamental form is given by:

$$(\nabla df)(X,Y) = \bar{\nabla}_{df(X)}\, df(Y) - df(\nabla_X Y).$$

The form $\nabla df$ is symmetric in $X$ and $Y$ and encodes the deviation of $f$ from being a totally geodesic map.

A smooth map $f: (M, g) \to (\bar{M}, \bar{g})$ between two Riemannian manifolds is said to be *totally geodesic* if its second fundamental form $\nabla df$ vanishes (see [2], [3]). In particular, any totally geodesic map has constant rank. Moreover, if $f$ is a totally geodesic map of constant rank $k$, with $0 < k < min\{m, n\}$, then the kernel of $df$ defines a smooth distribution with totally geodesic leaves, while its orthogonal complement is integrable with totally umbilical leaves (see [2], [3]).

Based on these observations, we have established the following result (see [6]):

**Theorem 2.1.** *Assume that $f$ is a smooth map of constant rank $k$, with $0 < k < n$, from a connected, compact (without boundary) Riemannian manifold $(M, g)$ of dimension $n > 2$ to a Riemannian manifold $(\bar{M}, \bar{g})$ of dimension $m \geq 2$. If $(M, g)$ has everywhere negative sectional curvature, then $f$ cannot be a totally geodesic map.*

The Dirichlet energy of a smooth map $f$ from a compact Riemannian manifold $(M, g)$ without boundary to a Riemannian manifold $(\bar{M}, \bar{g})$ is defined as (see [7]):

$$E_1(f) = \int_M \|df\|^2 \, dv_g,$$

where $dv := dv_g$ is the volume form of $M$ with respect to $g$. The differential $df$ is regarded as a section of the bundle of covariant 1-forms on $M$ with values in the pull-back of the tangent bundle of $\bar{M}$. Its magnitude is measured using the *Hilbert–Schmidt norm*, also known as the *Frobenius norm*, induced by the metrics $g$ and $\bar{g}$. That is, the squared norm of $df$ at each point of $M$ is computed relative to $g$ and $\bar{g}$.

There is a well-known formal correspondence between the energy density of a map and the norm of its differential (see [7]). In the context of Riemannian geometry, the energy density, denoted by $e(f)$, is defined as the trace of the pull-back metric $g^*$with respect to the metric $g$. For a map of constant rank, $e(f)$ quantifies the "stretching" of the map along the $k$ directions orthogonal to the fibers.

In particular, if $f$ is an isometric immersion, the pull-back metric $g^*$ coincides with $g$, and the energy density $e(f)$ equals n divided by two, where $n$ is the dimension of the domain manifold. More generally, in the Riemannian setting, the energy density is also expressed as one-half of the squared norm of the differential $df$. That is, the precise relationship can be written as $e(f)$ equals the squared norm of $df$.

Consequently, the squared norm of the differential $df$ can be expressed as the sum of the non-zero eigenvalues of the pull-back metric $g^*$. That is, the squared norm of $df$, which equals the energy density $e(f)$, is given by the sum of the eigenvalues $\lambda_1, \dots, \lambda_n$. The constant rank assumption guarantees that the number of non-zero terms in this sum remains fixed across the manifold. This prevents analytical singularities that typically arise when the rank of a map varies from point to point.

The Euler–Lagrange operator associated with the Dirichlet energy $E_1(f)$ is given by the trace of the second fundamental form $\nabla df$ with respect to the metric $g$, i.e.,

$$\tau_1(f) = trace_g\,(\nabla df).$$

This operator is called the *tension field* of the map $f$, denoted by $\tau_1(f)$.

A smooth map f from $(M, g)$ to $(\overline{M}, \bar{g})$ is said to be *harmonic* if its tension field $\tau_1(f)$ vanishes identically (see [7]).

Clearly, any totally geodesic map provides a trivial example of a harmonic map, since its tension field $\tau_1(f)$, given by the trace of the second fundamental form $\nabla df$, vanishes identically.

**Remark 2.1.** Further details on harmonic maps can be found in the monographs [8] and [9]. For instance, in the specific case where $f$ is a Riemannian submersion—a map of constant rank $k$ equal to the dimension of the target manifold $\overline{M}$, such that $df$ restricted to the horizontal distribution (the orthogonal complement of the kernel of $df$) is an isometry—the harmonicity of $f$ is closely related to the extrinsic geometry of its fibers. In particular, the mean curvature vector field of the fibers must vanish identically for $f$ to be harmonic.

We recall here the classical vanishing theorem for harmonic maps (see [5]; [7]).

**Theorem 2.2.** *Let $f$ be a harmonic map from a connected, compact (without boundary) Riemannian manifold $(M, g)$ with nonnegative Ricci curvature to a Riemannian manifold $(\overline{M}, \bar{g})$ with everywhere nonpositive sectional curvature. Then $f$ is necessarily totally geodesic and has constant rank. Moreover, if there exists at least one point in $M$ at which the Ricci curvature is strictly positive, then $f$ must be constant.*

**Example.** As an illustration of the Eells–Sampson theorem, consider a harmonic map $f$ from $(M, g)$ to $(\overline{M}, \bar{g})$ between two symmetric spaces of compact and noncompact types, respectively (see the definitions in [10]). The compact-type symmetric space is a connected, compact Riemannian manifold $(M, g)$ with nonnegative Ricci curvature. In the irreducible case, $(M, g)$ is an Einstein manifold with positive Ricci curvature. In contrast, the noncompact-type symmetric space is a simply connected, complete Riemannian manifold $(\bar{M}, \bar{g})$ with nonpositive sectional curvature.

It follows from Theorem 2 that any harmonic map from $(M, g)$ to $(\bar{M}, \bar{g})$ is totally geodesic. Moreover, if $(M, g)$ is irreducible, such a map must be constant.

The *bienergy* of a smooth map $f$ from a compact Riemannian manifold $(M, g)$ to a Riemannian manifold $(\bar{M}, \bar{g})$ is defined by (see [11]).

$E_2(f)$ is an integral over $M$ of the squared norm of the tension field $\tau_1(f)$ with respect to the volume element $dv_g$, i.e.,

$$E_2(f) := \int_M \|\tau_1(f)\|^2 \, dv_g.$$

The corresponding Euler–Lagrange operator associated with the bienergy functional $E_2(f)$ is given by

$$\tau_2(f) = \Delta_f \tau_1(f) - trace\, \bar{R}\big(df,\ \tau_1(f)\big)\, df,$$

where $\Delta_f$ denotes the rough Laplacian acting on smooth sections of the pull-back bundle $f^{-1}(T\bar{M})$, and $\bar{R}$ is the curvature tensor of the target manifold $\bar{M}$.

The vector field $\tau_2(f)$ is called the *bitension field* of the map $f$.

A smooth map $f$ is said to be biharmonic if its bitension field vanishes identically, that is, if $\tau_2(f) = 0$ (see [11]).

**Remark 2.2.** Further details on biharmonic mappings can be found in the monograph [12] and in the papers [11]; [13] and [14]. In particular, any harmonic map is a trivial example of a biharmonic map, since the vanishing of the tension field $\tau_1(f)$ immediately implies the vanishing of the bitension field $\tau_2(f)$.

Moreover, the following fundamental result holds (see [11]).

**Theorem 2.3.** *Let $f$ be a biharmonic map from a connected, compact Riemannian manifold $(M, g)$ into a Riemannian manifold $(\bar{M}, \bar{g})$ with everywhere non-positive sectional curvature. Then $f$ is harmonic.*

**Remark.** It should be emphasized that in the two main results on harmonic and biharmonic mappings, namely Theorems 2.2 and 2.3, the curvature assumption is imposed on the target manifold. In both cases, the non-positive sectional curvature of the image manifold plays a decisive role in establishing the corresponding energy inequalities and in deriving the associated rigidity conclusions.

## 3. A Dirichlet Energy–Bienergy Inequality for Maps into Manifolds of Non-Positive Sectional Curvature

The well-known Theorems 2.2 and 2.3 characterize harmonic and biharmonic maps from compact Riemannian manifolds into manifolds of non-positive sectional curvature. The following result provides a natural generalization of these classical theorems and constitutes the main theorem of the present paper.

**Theorem 3.1.**

*Let* $f: (M, g) \to (\bar{M}, \bar{g})$ *be a smooth map from a connected compact Riemannian manifold* $(M, g)$ *without boundary,* $dim\ M = n \geq 2$, *into a Riemannian manifold* $(\bar{M}, \bar{g})$, $dim\ \bar{M} = m \geq 2$, *whose sectional curvature is non-positive on the image* $f(M) \subset \bar{M}$. *Then the Dirichlet energy* $E_1(f)$ *and the bienergy* $E_2(f)$ *satisfy*

$$E_2(f) \geq Ric_{min}\, E_1(f),$$

*where* $Ric_{min}$ *denotes the global minimum of the Ricci curvature of* $(M, g)$.

*Moreover, equality holds if and only if* $f$ *is totally geodesic and has constant rank.*

**Proof.** On a compact Riemannian manifold without boundary, the following integral identity holds (see [15]):

$$\int_M (Q(f) + \| \nabla df \|^2 - \| \tau(f) \|^2)\, dv_g = 0. \tag{3.1}$$

Here $\tau(f)$ denotes the tension field, and $\| \nabla df \|^2$ is the squared norm of the second fundamental form. The scalar function $Q(f)$ is given by

$$Q(f) = \sum_{i=1}^{n} \lambda_i\ \mathrm{Ric}(e_i, e_i) - \sum_{i,j=1}^{n} \lambda_i \lambda_j\, \overline{\sec}(df(e_i), df(e_j)),$$

where $\{e_1, \dots, e_n\}$ is a local orthonormal frame on $M$, $\lambda_i$ are the eigenvalues of the pull-back metric $g^*$ with respect to $g$, and $\overline{\sec}$ denotes the sectional curvature of $(\bar{M}, \bar{g})$.

Since $E_2(f) = \int_M \| \tau(f) \|^2 \, dv_g$, identity (3.1) yields

$$E_2(f) = \int_M (Q(f) + \| \nabla df \|^2) \, dv_g \geq \int_M Q(f) \, d\, v_g.$$

Let $\mathrm{Ric}_{\min}$ be the minimum of the Ricci curvature of $(M, g)$. Then pointwise

$$\sum_{i=1}^{n} \lambda_i \ \mathrm{Ric}(e_i, e_i) \ \geq \ \mathrm{Ric}_{\min} \sum_{i=1}^{n} \lambda_i = \mathrm{Ric}_{\min} \ \| df \|^2.$$

If, in addition, $\overline{\sec} \leq 0$ on $f(M)$, then

$$-\sum_{i,j=1}^{n} \lambda_i \, \lambda_j \, \overline{\sec}(df(e_i), df(e_j)) \geq 0,$$

and hence

$$Q(f) \geq \mathrm{Ric}_{\min} \ \| df \|^2.$$

Integrating this inequality over $M$, we obtain

$$\int_M Q(f) \, d\, v_g \ \geq \ \mathrm{Ric}_{\min} \int_M \| df \|^2 \, dv_g.$$

By the definition of the Dirichlet energy,

$$E_1(f) = \int_M \| df \|^2 \, dv_g,$$

and therefore

$$\int_M Q(f) \, d\, v_g \ \geq \ \mathrm{Ric}_{\min} \, E_1(f).$$

On the other hand, from the integral identity (3.1)it follows that

$$E_2(f) = \int_M \| \tau(f) \|^2 \, dv_g = \int_M (Q(f) + \| \nabla df \|^2) \, dv_g \ \geq \ \int_M Q(f) \, d\, v_g.$$

Combining the above inequalities, we conclude that

$$E_2(f) \geq \mathrm{Ric}_{\min} \, E_1(f).$$

Assume now that equality holds:

$$E_2(f) = \mathrm{Ric}_{\min} \, E_1(f).$$

Then equality holds in all intermediate inequalities, and (3.1) reduces to

$$\int_M \Big(-\sum_{i,j=1}^{n} \lambda_i\,\lambda_j\,\overline{\mathrm{sec}}(df(e_i),df(e_j))\Big) + \Big(\sum_{i=1}^{n} \lambda_i\ \mathrm{Ric}(e_i,e_i) - \mathrm{Ric}_{\min} \parallel df \parallel^2\Big) +$$
$$\parallel \nabla df \parallel^2)\, dv_g = 0.$$

Each term in the integrand is nonnegative, hence all must vanish identically.
In particular, $\nabla df \equiv 0$, so $f$ is totally geodesic. Since $df$ is parallel and $M$ is connected, the rank of $df$ is constant on $M$. Moreover, $\mathrm{Ric}(e_i,e_i) = \mathrm{Ric}_{\min}$ for all $i$ with $\lambda_i \neq 0$, and $\overline{sec}\left(df(e_i),df\left(e_j\right)\right) = 0$ whenever $\lambda_i\lambda_j \neq 0$. Conversely, if $f$ is totally geodesic and has constant rank, then all inequalities above become equalities. □

**Remark 3.1.**

Equality forces the image of $df$ to lie in directions where the Ricci curvature of $(M, g)$ attains its minimum, while the sectional curvature of the target vanishes on all 2-planes contained in Im $df$. In addition, the inequality $E_2(f) \geq \mathrm{Ric}_{\min}\, E_1(f)$ can be interpreted as a spectral gap–type estimate for the tension field operator, in the spirit of Lichnerowicz's classical bound for the first eigenvalue of the Laplace–Beltrami operator. It shows that the $L^2$-norm of $\tau(f)$ is bounded below by $\mathrm{Ric}_{\min}$ times the $L^2$-norm of $df$, providing a coercivity property that geometrically controls the deviation from harmonicity and highlights the variational link between Dirichlet energy and bienergy.
On the other hand, the strict inequality

$$E_2(f) < Ric_{min} \cdot E_1(f)$$

is incompatible with the conclusion of Theorem 3.1. Indeed, if such a strict inequality were to hold, it would contradict the vanishing of the second fundamental form implied by equality in Theorem 3.1. We may therefore formulate the following vanishing theorem.

**Theorem 3.2.** *There exists no smooth map $f$ from a connected, compact Riemannian manifold $(M, g)$ of dimension $n \geq 2$ into a Riemannian manifold $(\overline{M}, \bar{g})$ of dimension $m \geq 2$ with non-positive sectional curvature at every point of $f(M)$ such that the Dirichlet energy $E_1(f)$ and the bienergy $E_2(f)$ of $f$ satisfy the strict inequality*

$$E_2(f) < Ric_{min}\, E_1(f)$$

*where $Ric_{min}$ denotes the global minimum of the Ricci curvature of $(M, g)$.*

**Remark 3.2.** If $f$ is a harmonic map, then the results of Theorems 4 and 5 immediately imply the classical vanishing theorem of Eells and Sampson for harmonic maps between compact Riemannian manifolds and manifolds of non-positive sectional curvature.

## 4. Special Case: Maps into Hadamard Manifolds

A *Hadamard manifold*, also known as a *Cartan–Hadamard manifold*, is a fundamental object in Riemannian geometry characterized by three properties: it is simply connected, complete, and has non-positive sectional curvature everywhere (see [L, p. 54]).

*Examples of Hadamard Manifolds*:

1. *Euclidean space* ($\mathbb{R}^n$): The simplest example, with zero sectional curvature.
2. Hyperbolic space ($\mathbb{H}^n$): A manifold with constant negative sectional curvature.
3. *Symmetric spaces of non-compact type:* Simply connected non-positively curved symmetric spaces are classical examples.

*Applications of Hadamard Manifolds*. Hadamard manifolds play an important role in modern computational and applied mathematics:

1. *Optimization:* They provide the standard setting for Riemannian Stochastic Gradient Descent and non-convex optimization problems (see, e.g., [17]).
2. *Machine Learning:* Used for dimensionality reduction and data analysis where the underlying data structures exhibit intrinsic non-positive curvature (see, e.g., [18]).
3. *Mathematical Physics:* Employed in the study of global spacetime structure and gravitational fields (see, e.g., [19]).

Next, we extend our investigation of the global geometry of Hadamard manifolds (see [20] and [21]) to smooth maps, focusing on their Dirichlet energy and bienergy.

The first proposition is an immediate corollary of Theorem 3.1.

**Corollary 4.1.** *Let $f$ be a smooth map from a connected, compact Riemannian manifold $(M, g)$ of dimension $n \geq 2$ into a Hadamard manifold $(\overline{M}, \bar{g})$. Then the Dirichlet energy $E_1(f)$ and the bienergy $E_2(f)$ of f satisfy the inequality*

$$E_2(f) \geq Ric_{min} \cdot E_1(f)$$

*where $Ric_{min}$ denotes the global minimum of the Ricci curvature of $(M, g)$. Moreover, equality holds if and only if $f$ is totally geodesic and has constant rank.*

Moreover, the second proposition is a direct consequence of Theorem 3.2.

**Corollary 4.2.** *There exists no smooth map $f$ from a connected, compact Riemannian manifold $(M, g)$ of dimension $n \geq 2$ into a Hadamard manifold $(\bar{M}, \bar{g})$ such that the Dirichlet energy $E_1(f)$ and the bienergy $E_2(f)$ of $f$ satisfy the strict inequality*

$$E_2(f) < Ric_{min} E_1(f),$$

*where $Ric_{min}$ denotes the global minimum of the Ricci curvature of $(M, g)$.*

## 5. Special Case: Totally Geodesic Maps

In differential geometry, projective maps between Riemannian or pseudo-Riemannian manifolds by definition preserve geodesics (see, for example, [M]). Studying them helps classify manifolds by their geodesic structure and explore symmetries beyond isometries, providing insight into the underlying geometry of space-time.

A map $\gamma: t \to \gamma(t)$ from an open interval $J \subset \mathbb{R}$ into a Riemannian manifold $(M, g)$ is said to be geodesic if it satisfies $\nabla_X X = \rho(t)X$, where $X = \acute{\gamma}$ Is tangent to $\gamma$. It is easily seen that the above equation can be reduced to $\nabla_X X = 0$ by a suitable change of the parameter $t$. In this case $\gamma$, is called geodesic with canonic (or affine) parameter $t$.

Suppose $f$ is a projective map of constant rank $k \leq min\,\{m, n\}$. Then following Yano (see [5]) we have

$$\nabla df = \theta \otimes d\,f + df \otimes \theta \tag{5.1}$$

where $\theta$ is a smooth 1-form on $(M, g)$.

In turn, we have proven the following statement in [15].

**Theorem 5.1**. *Let $f: (M, g) \to (\bar{M}, \bar{g})$ be a projective map of constant rank $k$ between closed oriented manifolds $(M, g)$ and $(\bar{M}, \bar{g})$ with $\overline{sec} \leq 0$ at each point of $f(M)$. If $Ric > 0$, then $f$ is a constant map. On the other hand, if $Ric \geq 0$, then $f$ is a map of constant rank $k \leq 1$.*

We next assume that an arbitrary geodesic in $(M, g)$ is mapped by $f$ into a geodesic in $(\bar{M}, \bar{g})$ such that the affine parameter is preserved. Such a mapping $f$ is said to be totally geodesic or affine (see [5] and [2]). In order for a mapping $f:(M,g)\to(\overline{M},\bar{g})$ to be totally geodesic (or affine), it is necessary and sufficient that $\nabla df=0$ (see [5] and [2]). Our main result of this paragraph is the following.

**Theorem 5.2**. *Let $f:(M,g)\to(\overline{M},\bar{g})$ be a projective map of constant rank $k$ between a connected, compact Riemannian manifold $(M,g)$ of dimension $n\geq 2$ and a Riemannian manifold $(\overline{M},\bar{g})$ of dimension $m\geq 2$ with non-positive sectional curvature at each point of $f(M)$. Then the Dirichlet energy $E_1(f)$ and the bienergy $E_2(f)$ of $f$ satisfy the inequality*

$$E_2(f)\geq 2\,Ric_{min}\cdot E_1(f),$$

*where $Ric_{min}$ denotes the global minimum of the Ricci curvature of $(M,g)$. Moreover, equality holds if $f$ is a constant map.*

**Proof**. Let $f:(M,g)\to(\overline{M},\bar{g})$ be a projective map of constant rank from a compact Riemannian manifold $(M,g)$ of dimension $n\geq 2$ to a Riemannian manifold $(\overline{M},\bar{g})$ of dimension $m\geq 2$. In this case (5.1) holds. Then the tension field of the map $f$ has the form

$$\tau_1(f)=trace_g\,(\nabla df)=2\,df(\xi),$$

where $\xi=\theta^{\#}$ and, therefore, we have

$$E_2(f):=\int_M\|\tau_1(f)\|^2dv_g=4\int_M\|df(\xi)\|^2\,dv_g$$

In this case from (3.1) we obtain

$$E_2(f)=\ 2\int_M(Q(f)+2\|\theta\|^2\|df\|^2)\,dv_g. \tag{5.2}$$

From (5.2) we deduce

$$E_2(f)\geq\ 2\int_M Q(f)\,dv_g, \tag{5.3}$$

where the inequality

$$Q(f)(x)\geq-\sum_{i,j=1}^{k}\lambda_i\lambda_j\,\overline{sec}\left(df(e_i),df(e_j)\right)+Ric_{min}\left(\sum_{i=1}^{k}\lambda_i\right).$$

holds at each point $x\in M$.

If, in addition, $(\overline{M},\bar{g})$ is a Riemannian manifold with non-positive sectional curvature, then $Q(f)(x)\geq Ric_{min}\cdot\|df\|^2$. Then from (5.3) we deduce the following inequality:

$$E_2(f)\geq\ 2\,Ric_{max}\cdot E_1(f).$$

In particular, if we suppose that

$$E_2(f) = 2\, Ric_{max} \cdot E_1(f) =$$
$$= 2\, Ric_{max} \int_M \|df\|^2\, dv_g = 2\, Ric_{max} \int_M \left(\sum_{i=1}^{k} \lambda_i\right) dv_g,$$

then we can rewrite the integral identity (5.2) in the form

$$\int_M \left(-\sum_{i,j=1}^{n} \lambda_i \lambda_j\, \overline{sec}\left(df(e_i), df(e_j)\right) + \left(\sum_{i=1}^{n} \lambda_i Ric(e_i, e_i) - Ric_{min}\left(\sum_{i=1}^{k} \lambda_i\right)\right) + 2\|\theta\|^2 \|df\|^2\right) dv_g = 0.$$

Each of the three terms in the integrand is nonnegative. Hence, the vanishing of the integral implies that each term must vanish identically. In particular, we obtain that $\theta = 0$ and $df = 0$. Then in the case of a connected Riemannian manifold $(M, g)$, the map $f$ is constant. ■

## 6. Application: The functional $\mathcal{F}_{\lambda,\mu}(f) = \mu E_1(f) + \lambda E_2(f)$

The inequality obtained in Theorem 3.1 admits a natural interpretation within the variational framework related to sesqui-harmonic maps. To this end, consider the functional

$$\mathcal{F}_{\lambda,\mu}(f) = \mu E_1(f) + \lambda E_2(f),$$

where $\lambda, \mu \in \mathbb{R}$ and $f: (M, g) \to (\bar{M}, \bar{g})$is a smooth map between Riemannian manifolds. Its Euler–Lagrange equation is given by

$$\mu\, \tau_1(f) + \lambda\, \tau_2(f) = 0,$$

and its critical points are called sesqui-harmonic (or interpolating) maps.

**Remark 6.1.** The functional $F_{\lambda,\mu}(f)$ is closely related to earlier developments in higher-order variational calculus (see, e.g., [23]). In particular, the functional combining the Dirichlet energy and the bienergy arises naturally in the theory of sesqui-harmonic maps introduced by Branding [24]. Related variational problems have also been studied in the context of biminimal immersions (see [25]).

The inequality established in Theorem 3.1 of the present paper provides a direct estimate for this functional.

**Theorem 6.1.** *Let* $f:(M,g)\to(\bar{M},\bar{g})$ *be a smooth map from a connected compact Riemannian manifold* $(M,g)$ *into a Riemannian manifold* $(\bar{M},\bar{g})$*whose sectional curvature is non-positive on the image* $f(M)\subset\bar{M}$. *Then the functional* $F_{\lambda,\mu}(f)=\mu E_1(f)+\lambda E_2(f)$ *satisfies:*

(i) *If* $\lambda\geq 0$, *then*

$$F_{\lambda,\mu}(f)\geq(\mu+\lambda\,Ric_{min})\,E_1(f).$$

*In particular, if* $\mu+\lambda\,Ric_{min}\geq 0$, *then* $F_{\lambda,\mu}(f)\geq 0$.

(ii) *If* $\lambda<0$, *then*

$$F_{\lambda,\mu}(f)\leq(\mu+\lambda\,Ric_{min})\,E_1(f).$$

*In particular, if* $\mu+\lambda\,Ric_{min}\leq 0$, *then* $F_{\lambda,\mu}(f)\leq 0$.

**Proof.** By Theorem 3.1, we have $E_2(f)\geq\mathrm{Ric}_{\min}\,E_1(f)$.

Firstly, if $\lambda\geq 0$, multiplication preserves the inequality, giving $\lambda E_2(f)\geq\lambda\,\mathrm{Ric}_{\min}\,E_1(f)$, and hence

$$F_{\lambda,\mu}(f)\geq(\mu+\lambda\,\mathrm{Ric}_{\min})\,E_1(f).$$

Therefore, if $\mu+\lambda\,\mathrm{Ric}_{\min}\geq 0$, then $F_{\lambda,\mu}(f)\geq 0$.

Secondary, if $\lambda<0$, multiplication reverses the inequality, yielding

$$F_{\lambda,\mu}(f)\leq(\mu+\lambda\,\mathrm{Ric}_{\min})\,E_1(f).$$

Therefore, if $\mu+\lambda\,\mathrm{Ric}_{\min}\leq 0$, then $F_{\lambda,\mu}(f)\leq 0$. □

**Corollary 6.1** (Rigidity). *Let* $f:(M,g)\to(\bar{M},\bar{g})$ *be as above. Assume that* $\lambda\neq 0$ *and that equality holds in Theorem 6.1, i.e.,*

$$F_{\lambda,\mu}(f)=(\mu+\lambda\,Ric_{min})\,E_1(f).$$

*Then* $f$ *is totally geodesic and has constant rank.*

**Proof.** Equality implies $E_2(f)=\mathrm{Ric}_{\min}\,E_1(f)$. Then the integral identity (3.1) takes the form

$$\int_M \left(-\sum_{i,j=1}^{n} \lambda_i \lambda_j\, \overline{sec}(df(e_i), df(e_j)) + \left(\sum_{i=1}^{n} \lambda_i \ \mathrm{Ric}(e_i, e_i) - \mathrm{Ric}_{\min} \parallel df \parallel^2\right) + \right.$$

$$\left. \parallel \nabla df \parallel^2 \right) dv_g = 0.$$

Each term in the integrand is nonnegative. Therefore, each of them must vanish identically on $M$. In particular, the vanishing of the third term yields $\nabla df \equiv 0,$ so $f$is totally geodesic. Since $df$is parallel and $M$is connected, the rank of $df$is constant on $M$. Furthermore, the vanishing of the second term implies that, at each point of $M$, one has $\mathrm{Ric}(e_i, e_i) = \mathrm{Ric}_{\min}$ for every direction $e_i$corresponding to a nonzero singular value $\lambda_i$. Finally, the vanishing of the first term shows that

$$\overline{\mathrm{sec}}(df(e_i), df(e_j)) = 0$$

whenever $\lambda_i \lambda_j \neq 0$, that is, the sectional curvature of the target vanishes on every 2-plane contained in the image of $df$. □

**Remark 6.2.** Theorem 6.1 shows that the Ricci curvature of the domain provides a natural geometric control of the functional $F_{\lambda,\mu}(f)$i n terms of the Dirichlet energy. In particular, for $\lambda \geq 0$and $\mu \geq 0$, the functional $F_{\lambda,\mu}(f)$ admits a nonnegative lower bound proportional to $E_1(f)$, reflecting a coercivity-type property under curvature constraints.

## 7. Conclusion

In this paper, we established a global inequality relating the bienergy and the Dirichlet energy of smooth maps between Riemannian manifolds under natural curvature assumptions. Specifically, if $(M, g)$is a compact Riemannian manifold whose Ricci curvature admits a global lower bound $\mathrm{Ric}_{\min}$, and if the target manifold has non-positive sectional curvature along the image of a smooth map $f$, then

$$E_2(f) \geq Ric_{min}\, E_1(f).$$

This inequality provides a quantitative relationship between the first-order and second-order energy functionals and shows that the Ricci curvature of the domain controls the bienergy from below in terms of the Dirichlet energy. In particular, this estimate has the character of a spectral gap–type inequality, expressing a coercivity property of the tension field in the $L^2$-sense. As a consequence, when $Ric_{min} > 0$, the bienergy dominates the energy in a strict sense, reinforcing rigidity phenomena and providing a geometric obstruction to the existence of nontrivial maps with relatively small bienergy. Thus, the obtained result clarifies the geometric and variational relationship between harmonic and biharmonic structures and provides a useful analytical tool for further studies of biharmonic maps and related higher-order variational problems.

**Sergey Stepanov**
Department of Mathematics and Data Analysis, Finance University, 49-55, Leningradsky Prospect, 125468 Moscow, Russia;
Department of Mathematics, Russian Institute for Scientific and Technical Information of the Russian Academy of Sciences, 20, Usievicha street, 125190 Moscow, Russia;
s.e.stepanov@mail.ru

**Irina Tsyganok**
Department of Mathematics and Data Analysis, Finance University, 49-55, Leningradsky Prospect, 125468 Moscow, Russia;
i.i.tsyganok@mail.ru